\documentclass[a4paper, 11pt]{article}
\usepackage[english]{babel}
\usepackage{latexsym}
\usepackage{amsfonts, amsthm}
\usepackage{amssymb}
\usepackage{amscd}
\usepackage{amsmath}
\usepackage{graphics}
\usepackage{eepic}
\usepackage{epsfig}
\usepackage[noxcolor]{pstricks}
\newtheorem{theorem}{Theorem}

\newtheorem{remark}{Remark\/}

\newcommand{\C}{{\mathbb C}}
\newcommand{\R}{{\mathbb R}}

\newcommand{\Z}{{\mathbb Z}}

\newcommand{\<}{\backslash}

\begin{document}

\title{\textbf{Poles of the topological zeta function for plane curves and Newton polyhedra}}
\author{Ann Lemahieu and Lise Van Proeyen \footnote{Ann Lemahieu, Lise Van Proeyen; K.U.Leuven, Departement Wiskunde,
Celestij\-nenlaan 200B, B-3001 Leuven, Belgium, email:
lemahieu@mathematik.uni-kl.de, lise.vanproeyen@wis.kuleuven.be. The
research was partially supported by the Fund of Scientific Research
- Flanders (G.0318.06) and MEC PN I+D+I MTM2007-64704.}
\date{}}
\maketitle {\footnotesize \emph{\noindent \textbf{Abstract.---} The
local topological zeta function is a rational function associated to
a germ of a complex holomorphic function. This function can be
computed from an embedded resolution of singularities of the germ.
For nondegenerate functions it is also possible to compute it from
the Newton polyhedron. Both ways give rise to a set of candidate
poles of the topological zeta function, containing all poles. \\
\indent For plane curves, Veys showed how to filter the actual poles
out of the candidate poles induced by the resolution graph. In this
note we show how to determine from the Newton polyhedron of a
nondegenerate plane curve which candidate poles are actual poles.}}
\\ \\
 ${}$ \begin{center}
\textsc{1. Introduction}
\end{center} ${}$\\
\textbf{1.1 The local topological zeta function.--\quad} In 1992
Denef and Loeser introduced a new zeta function which they called
the topological zeta function because of the topological
Euler--Poincar\'e characteristic turning up in it. Let $f: (\C^n,0)
\rightarrow (\C,0)$ be the germ of a holomorphic function and let
$\pi : X \rightarrow \C^n$ be an embedded resolution of $f^{-1} \{ 0
\}$. We denote by $E_i, i \in S$, the irreducible components of
$\pi^{-1}(f^{-1}\{ 0 \})$, and by $N_i$ and $\nu_i - 1$ the
multiplicities of $E_i$ in the divisor on $X$ of $f \circ \pi$ and
$\pi^\ast (dx_1 \wedge \ldots \wedge dx_n)$, respectively. For $I
\subset S$ we denote also $E_I := \cap_{i \in I} E_i$ and $E^\circ_I
:= E_I \setminus (\cap_{j \notin I} E_j)$. Further we write
$\chi(\cdot)$ for the topological Euler--Poincar\'e characteristic.

The \emph{local topological zeta function associated to $f$} is the
rational function in one complex variable $Z_{top,f} (s) := \sum_{I
\subset S} \chi (E^\circ_I \cap h^{-1}\{ 0 \}) \prod_{i \in I}
\frac{1}{N_i s+ \nu_i}$. Denef and Loeser proved in \cite{DL-zeta
fie onafh van res} that these definitions are independent of the
choice of the resolution.

The poles of the local topological zeta function are part of the set
$\{-\nu_i/N_i \mid i \in S\}$; therefore this set is called a set of
candidate poles. Various conjectures, such as the monodromy
conjecture and the holomorphy conjecture, relate the poles to the
eigenvalues of the local monodromy of $f$ (see for example
\cite{Denef Bourbaki rapport}). A very remarkable fact is that most
of the candidate poles are cancelled in the topological zeta
function. In general, it is not known how to see whether a candidate
pole is a pole or not. Only for plane curves there exists a complete
criterion. In \cite{veyspolen} Veys showed:

\begin{theorem} Let $f \in
\mathbb{C}[x,y]$ be a non-constant polynomial satisfying $f(0)=0,$
and let $\pi:X \to \mathbb{C}^2$ be the minimal embedded resolution
of $f^{-1}\{0\}$ in a neighbourhood of 0. Then $s_0$ is a pole of
$Z_{top,f}(s)$ if and only if $s_0 = -\frac{\nu_i}{N_i}$ for some
exceptional curve $E_i$ intersecting at least three times other
components or $s_0 = -\frac1{N_i}$ for some irreducible component
$E_i$ of the strict transform of $f=0.$
\end{theorem}

\noindent \textbf{1.2 The topological zeta function out of the
Newton polyhedron.--} Let $f \in \C[x_1, \ldots , x_n]$ be a
non-constant polynomial satisfying $f(0)=0.$ We write $f = \sum_{k
\in \Z_{\geq 0}^n} a_k x^k,$ where $k=(k_1, \ldots, k_n)$ and $x^k =
x^{k_1} \cdot \ldots \cdot x_n^{k_n}.$ The \emph{support of $f$} is
$\textnormal{supp}\, f := \{k \in \Z_{\geq 0}^n \, | \, a_k \neq 0
\}.$ The \emph{Newton polyhedron $\Gamma_0$ of $f$ at the origin} is
the convex hull in $\mathbb{R}_{\geq 0}^n$ of $\bigcup_{k \in
\textnormal{{supp}}\, f} k + \mathbb{R}_{\geq 0}^n.$ A \emph{face}
of the Newton polyhedron is the intersection of $\Gamma_0$ with a
supporting hyperplane. A \emph{facet} is a face of dimension $n-1.$
A polynomial $f(x_1, \ldots , x_n)$ is called \emph{nondegenerate
with respect to its Newton polyhedron $\Gamma_0$} if
 for every compact face $\tau$ of $\Gamma_0$
the polynomials $f_\tau := \sum_{k \in \tau} a_k x^k$ and $\partial
f_\tau /
\partial x_i, 1 \leq i \leq n,$ have no common zeroes in $(\C\<
\{0\})^n.$ 

If $a_1, \ldots , a_r \in \R^n \< \{0\},$ then
$\textnormal{cone}(a_1, \ldots , a_r):= \{\sum_{i=1}^r \lambda_i a_i
\, | \, \lambda_i \in \R, \lambda_i > 0 \}$. If $\Delta$ can be
written as $\textnormal{cone}(a_1, \ldots , a_r)$, with $a_1 ,
\ldots , a_r \in \Z^n \< \{0\}$ linearly independent over $\R$, then
$\Delta$ is called a \emph{simplicial cone.} For a simplicial cone
$\Delta$ spanned by the primitive and linearly independent vectors
$a_1 , \ldots , a_r \in \Z^n$, the \emph{multiplicity} of $\Delta,$
denoted by $\textnormal{mult}(\Delta)$, is the index of the lattice
$\Z a_1 + \ldots + \Z a_r$ in the group of the points with integral
coordinates of the vector space generated by $a_1 , \ldots, a_r.$ 
Then mult$(\Delta)$ is equal to the greatest common divisor of the
determinants of the $(r \times r)$-matrices obtained by omitting
columns from the matrix $A$ with rows $a_1, \ldots, a_r.$ 

Let $\Gamma_0$ be a Newton polyhedron in $\R^n$. For $a= (a_1,
\ldots , a_n) \in \R_{\geq 0}^n$ we put $N(a) := \textnormal{inf}_{x
\in \Gamma_0} a \cdot x, \nu (a) := \sum_{i=1}^n a_i$ and $F(a) :=
\{x \in \Gamma_0 \, | \, a \cdot x = N(a)\}.$ All $F(a), a \neq 0,$
are faces of $\Gamma_0.$ To a face $\tau$ of $\Gamma_0$ one
associates a dual cone $\tau^\circ \subset \R^n,$ defined as the
closure in $\R^n$ of $\{a \in \R_{\geq 0}^n \, | \, F(a) = \tau \}.$
This is a cone of dimension $n - \textnormal{dim } \tau$ with vertex
in the origin. For a facet $\tau$, one has $\tau^\circ = a \R_{\geq
0}$ for some primitive $a \in \Z_{\geq 0}^n,$ and then the equation
of the hyperplane through $\tau$ is $a \cdot x = N(a).$ We also use
the notation $N(\tau)$ and $\nu(\tau)$, meaning respectively $N(a)$
and $\nu(a)$ for this associated $a \in \Z_{\geq 0}^n.$ The set
$\{\tau^\circ \mid \tau \textnormal{ face of } \Gamma_0\}$ defines a
subdivision of $\R^n_{\geq 0}$ and is called the \emph{normal fan}
to $\Gamma_0$. In 1976 Varchenko proved in \cite{Varchenko} that the
map from the toric variety corresponding to a regular subdivision of
the normal fan to $\C^n$ is an embedded resolution for all
polynomials having $\Gamma_0$ as Newton polyhedron in the origin and
that are nondegenerate with respect to $\Gamma_0$. Denef and Loeser
used this to provide a formula for the local topological zeta
function out of the Newton polyhedron.

Suppose $\Delta = \R_{\geq 0}a_1 + \cdots + \R_{\geq 0}a_r,$ with
$a_1, \ldots , a_r \in \Z_{\geq 0}^n$ linearly independent and
primitive. They define
$J_\Delta(s):=\frac{\textnormal{mult}(\Delta)}{\prod_{i=1}^r(N(a_i)s
+ \nu(a_i))}$ and to an arbitrary face $\tau$ of $\Gamma_0$ they
associate the rational function $J_{\tau}(s) := \sum_{i=1}^k
J_{\Delta_i}(s)$, with $\tau^\circ = \cup_{i=1}^k \Delta_i$ a
decomposition of $\tau^\circ$ into simplicial cones $\Delta_i$ of
dimension $\ell = \textnormal{dim }\tau^\circ$ satisfying
$\textnormal{dim}(\Delta_i
\cap \Delta_j) < \ell$ if $i\neq j$.  

\begin{theorem} \label{stelling topzeta newton} \emph{\cite[Th\'eor\`eme 5.3]{DL-zeta fie onafh van res}}
If $f$ is nondegenerate with respect to $\Gamma_0,$ then the local
topological zeta function is equal to \emph{
$$Z_{top,f}(s) = \sum_{\tau \textnormal{ {vertex of }}\Gamma_0} J_\tau(s)
+ \frac{s}{s+1} \sum_{\substack{\tau \textnormal{ compact}\\
\textnormal{face of }\Gamma_0,\\
\textnormal{dim }\tau \geq 1}} (-1)^{\textnormal{dim
}\tau}(\textnormal{dim }\tau)! \textnormal{Vol}(\tau) J_\tau(s).$$}
\end{theorem}
\noindent For a face $\tau$ of dimension 0, $\textnormal{Vol}(\tau)
:= 1.$ For every other compact face $\textnormal{Vol}(\tau)$ is the
\emph{volume} of $\tau$ for the volume form $\omega_\tau.$ This is a
volume form on Aff$(\tau),$ the affine space spanned by $\tau,$ such
that the parallelepiped spanned by a lattice-basis of $\Z^n \cap
\textnormal{Aff}(\tau)$ has volume 1. The product $(\textnormal{dim
}\tau)! \textnormal{Vol}(\tau)$ is also called the \emph{normalized
volume of }$\tau.$ If $\tau$ is a simplicial facet, this normalized
volume is equal to the multiplicity of the cone spanned by the
vertices divided by $N(\tau)$.
\\ \\
Theorem \ref{stelling topzeta newton} yields another set of
candidate poles (containing all poles) of the local topological zeta
function, namely -1 together with the rational numbers
$-\nu(\tau)/N(\tau)$ for $\tau$ a facet of $\Gamma_0$. We will say
that such a facet \emph{contributes} the candidate pole. In the
following section we will give a criterion for a candidate pole of
this set to be a pole of the local topological zeta function.

\begin{remark}\label{remarkgrafisch}
\emph{There is also a graphical way to determine the candidate pole
contributed by a facet $\tau$. If $(r, \ldots , r)$ is the
intersection point of the diagonal of the first quadrant with the
affine hyperplane containing $\tau$, then the candidate pole
$-\nu(\tau)/N(\tau)$ is equal to $-1/r$.}
\end{remark}

 ${}$ \begin{center}
\textsc{2. Description of the poles in terms of the Newton
polyhedron}
\end{center}
${}$\\ We will say that a facet of a $2$-dimensional Newton
polyhedron is a \emph{$B_1$-facet} with respect to the variable $x$
(resp. to the variable $y$) if it has one vertex in the coordinate
hyperplane $x=0$ (resp. $y=0$) and one vertex at distance one of
this hyperplane.
\begin{theorem}
Let $f$ be a complex polynomial in two variables. Suppose that $f$
is nondegenerate with respect to its Newton polyhedron $\Gamma_0$.
Then for a candidate pole $s_0 \neq -1$ contributed by some facet of
$\Gamma_0$ it holds: $s_0$ is a pole of $Z_{top,f}$ if and only if
$s_0$ is contributed by a facet of $\Gamma_0$ that is no
$B_1$-facet.
\end{theorem}

\noindent \emph{Proof.} Suppose first that $s_0$ is a candidate pole
of order $2$. Let $\sigma$ and $\tau$ be facets such that
$s_0=\nu(\tau)/N(\tau)=\nu(\sigma)/N(\sigma)\neq 1$, having exactly
one point in common. According to Remark \ref{remarkgrafisch}, the
picture should
be as in Figure 1.
\begin{figure}[h]
\begin{center}
\resizebox{1.8in}{!}{\includegraphics{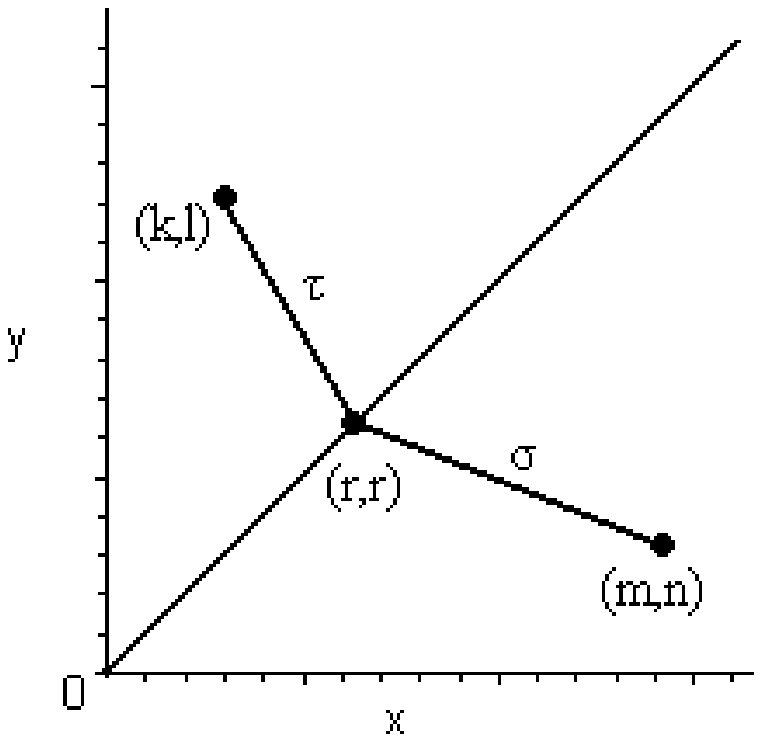}}
\resizebox{1.8in}{!}{\includegraphics{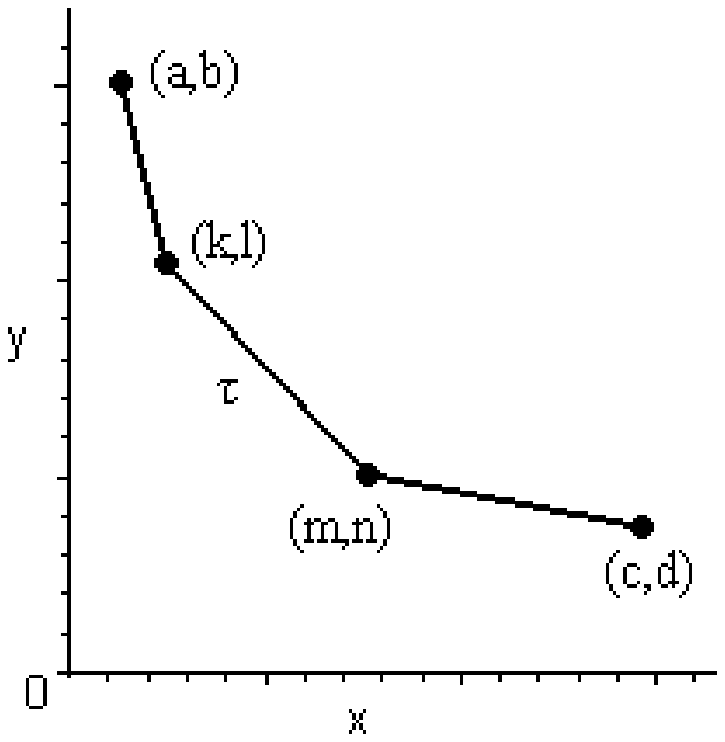}}
\resizebox{1.8in}{!}{\includegraphics{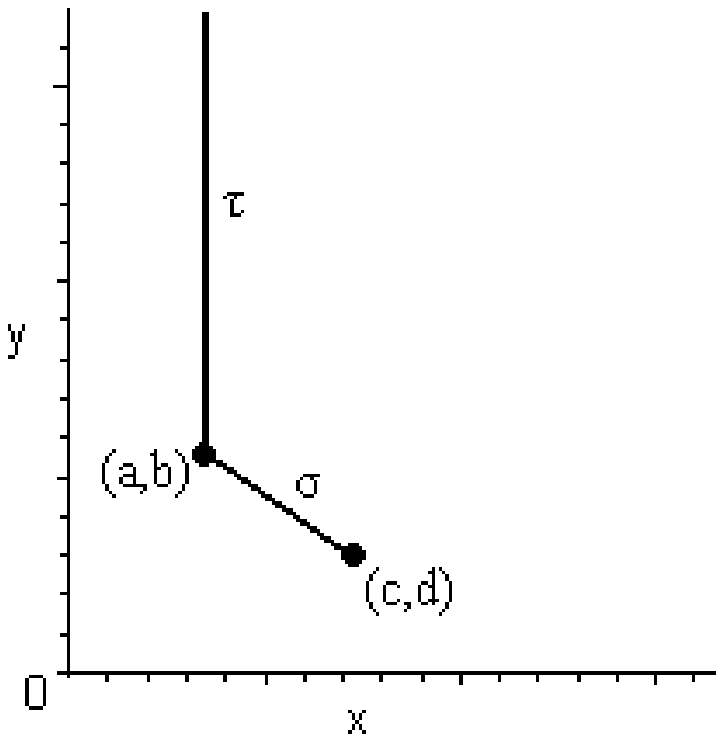}}\\
Figure 1 \qquad \qquad \qquad \qquad Figure 2 \qquad \qquad \qquad
\qquad Figure 3
\end{center}
\vspace*{-0.5cm}
\end{figure}
Note that $\tau$ or $\sigma$ might not be compact. The facets $\tau$
and $\sigma$ are no $B_1$-facets, since otherwise $r$ should be
equal to 1 and $s_0 = -1$. When looking at the formula for the
topological zeta function given in Theorem \ref{stelling topzeta
newton}, it is obvious that a candidate pole of order two is also a
pole of order two.


\vspace{\baselineskip}

Suppose now that $s_0$ is a candidate pole of order 1 and suppose
first that $s_0$ is contributed by a compact facet $\tau$. Then we
fix notation as in Figure 2. Note that $k\neq l$ and $m\neq n$.
We sum the contributions to the local topological zeta
function coming from the segment $\tau$ and the points with
coordinates $(k,l)$ and $(m,n)$ (see Theorem \ref{stelling topzeta
newton}). If $g$ is the normalized volume of $\tau$ (which implies
$g=\textnormal{gcd}(m-k,l-n)$), the contribution of $\tau$ is
\[-\frac{sg^2}{(s+1)((lm-kn)s+l-n+m-k)},\] the contribution of the
point with coordinates $(k,l)$ is given by
\[\frac{(b-l)(m-k)-(l-n)(k-a)}{((lm-kn)s+l-n+m-k)((bk-al)s+b-l+k-a)}\]
and the contribution of $(m,n)$ is
\[\frac{(l-n)(c-m)-(n-d)(m-k)}{((lm-kn)s+l-n+m-k)((nc-md)s+n-d+c-m)}.\]
Once summed these 3 contributions, we compute the residue Res at
$s=-(l-n+m-k)/(lm-kn)$: 

{\footnotesize \[{\mbox{Res } = \frac{(l-n +m
-k)((ml-nk)(ml-nk+k-m+n-l)+g^2(n-m)(k-l))}{
 (lm-kn)(n-m)(k-l)(ml-nk+k - m+n  -l)}}.\]}
If $\tau$ is a $B_1$-facet, then $k=0, m=1, g=1$ or $n=0, l=1, g=1.$
It is easy to calculate that then Res $=0$. From now on, we suppose
$\tau$ is no $B_1$-facet. First we note that $ml-nk+k-m+n-l>0$. This
follows from the fact that a candidate pole in dimension two is
always bigger than or equal to $-1$. Since we suppose $s_0 \neq -1$,
this leads to $\frac{l-n+m-k}{lm-kn} < 1$.


If $\tau$ is intersecting the diagonal of the first quadrant, then
there can clearly not exist another facet $\tau$ yielding this
candidate pole.
To show that the candidate pole $s_0=-\nu(\tau)/N(\tau)$ is a pole
of $Z_{top,f}$, it is thus sufficient to prove that Res $\neq 0$. As
in this situation $k<l, k<m, n<m, n<l$, it follows that Res $=0$ if
and only the factor
\[F:=(ml-nk)(ml-nk+k-m+n-l)+g^2(n-m)(k-l)=0.\]
We have $ml-nk > 0, k-l<0, n-m<0$ and $ml-nk+k-m+n-l >0$. This leads
to $F>0$ and Res $\neq 0$.


Suppose now that $\tau$ is lying above the diagonal in the first
quadrant. Then we have $k<m<n<l$. We prove now that Res $<0$. It is
easy to see that this is equivalent to $F>0$. 
We observe that
\begin{eqnarray*}
& & (ml-nk)(ml-nk+k-m+n-l)-(l-k)(n-m)(m-k)^2 \\
& = & lm(m-k)(m-k-1) + k(n-m)(m-k)^2+kn(m-k)\\
& & + (l-n)(lm(m-1)-nk^2)+kn(l-n).
\end{eqnarray*}
All terms in this summand are greater or equal than $0$ and we find
that the total expression is equal to $0$ exactly when $k=0$ and
$m=1$ or equivalently, when $\tau$ is a $B_1$-facet. If $\tau$ is
not a $B_1$-facet, then we find $F>0$ because $g^2 \leq (m-k)^2$.


If $\tau$ is a facet lying under the diagonal, then one can permute
$n$ with $k$ and $m$ with $l$. We get the same conclusion: the
residue is strictly negative unless the facet is a $B_1$-facet.

Now we are left with the case that $s_0$ is a candidate pole of
order 1, contributed by a facet $\tau$ that is not compact. Suppose
$\textnormal{Aff}(\tau) \leftrightarrow x = a$ and $(a,b),$ $(c,d)$
and $\sigma$ are as in Figure 3. (The line segment $\sigma$ might
not be compact, in that case $d=b$.)
Then the candidate pole $s_0 = -1/a$ only turns up in the term of
the local topological zeta function originating from the point
$(a,b)$. As $\tau$ induces a candidate pole of order 1, one has that
$a \neq b$. This term is equal to
$$\frac{c-a}{(as+1)((bc-ad)s+b-d+c-a)}.$$ The residue in $s_0$ is
$\frac{1}{(a-b)}.$ As in the case of a compact facet, when $\tau$
intersects the diagonal of the first quadrant, then $s_0$ is only
contributed by $\tau$ and thus $s_0$ is a pole. If $\tau$ is lying
completely above the diagonal, then $\frac{1}{(a-b)} < 0$. Also for
non-compact facets parallel with the $x$-axis one gets the same
conclusion.

Hence, contributions coming from different facets (compact or not
compact) do not cancel each other. This ends the proof. \qed

\begin{remark}
\emph{Notice that the computed residues do not depend on the
neighbour segments of $\tau$, although they are used in the
computation of the residue.}
\end{remark}

\footnotesize{

}


\begin{thebibliography}{AKMW}

\bibitem[De]{Denef Bourbaki rapport} J. Denef, \emph{Report on Igusa's local zeta
function,} S\'em. Bourbaki 741, Ast\'erisque \textbf{201/202/203}
(1991), 359-386.

\bibitem[DL]{DL-zeta fie onafh van res} J. Denef and F. Loeser,
\emph{Caract\'eristiques d' Euler-Poincar\'e, fonctions z\^eta
locales, et modifications analytiques,} J. Amer. Math. Soc.
\textbf{5}  (1992), 705-720.

\bibitem[Var]{Varchenko} A.N. Varchenko, \emph{Zeta-Function of Monodromy and Newton's
Diagram,} Inventiones math. \textbf{37} (1976), 253-262.

\bibitem[Ve]{veyspolen} W. Veys, \emph{Determination of the poles of the topological zeta
function for curves}, Manuscripta Math. \textbf{87} (1995), 435-448.

\end{thebibliography}
\end{document}